\documentclass{article}
\usepackage{latexsym}
\usepackage{epsfig}
\usepackage{amsmath}
\usepackage{amssymb}
\usepackage{amsthm}
\newtheorem{theorem}{Theorem}[section]
\newtheorem{corollary}[theorem]{Corollary}

\newtheorem{lemma}[theorem]{Lemma}

\theoremstyle{remark}
\newtheorem{remark}[theorem]{Remark}

\theoremstyle{definition}
\newtheorem{definition}[theorem]{Definition}

\DeclareMathOperator\diag{diag}

\begin{document}

\title{Natural parameterizations of \\ closed projective plane curves}

\author{Roland Hildebrand \thanks{%
Univ.~Grenoble Alpes, CNRS, Grenoble INP, LJK, 38000 Grenoble, France
({\tt roland.hildebrand@univ-grenoble-alpes.fr}).}}

\maketitle

\begin{abstract}
A natural parametrization of smooth projective plane curves which tolerates the presence of sextactic points is the Forsyth-Laguerre parametrization. On a closed projective plane curve, which necessarily contains sextactic points, this parametrization is, however, in general not periodic. We show that by the introduction of an additional scalar parameter $\alpha \leq \frac12$ one can define a projectively invariant $2\pi$-periodic global parametrization on every simple closed convex sufficiently smooth projective plane curve without inflection points. For non-quadratic curves this parametrization, which we call balanced, is unique up to a shift of the parameter. The curve is an ellipse if and only if $\alpha = \frac12$, and the value of $\alpha$ is a global projective invariant of the curve. The parametrization is equivariant with respect to duality.
\end{abstract}

Keywords: projective plane curve, Forsyth-Laguerre parametrization, global invariant

MSC 2010: 53A20, 52A10

\section{Parameterizations of projective plane curves} \label{sec:intro}

Projective plane curves have been intensely studied in the second half of the 19-th and the beginning of the 20-th century and are a classical subject of differential geometry. In this paper we consider periodic parameterizations of closed projective plane curves. The well-known natural local parameterizations cannot in general be extended to the whole curve. We show that under some non-degeneracy assumptions there nevertheless exists a natural periodic global parametrization. On non-quadratic curves it gives rise to a projectively invariant metric on the curve.

\medskip

The most natural way to represent curves in the real projective plane is by projective images of vector-valued solutions of third-order linear differential equations. This representation has already been studied in the 19-th century by Halphen, Forsyth, Laguerre, and others. For a detailed account see \cite{Wilczynski} or \cite{Cartan37}, for a more modern exposition see \cite{OvsienkoTabachnikov}.

Let $\gamma$ be a regularly parameterized (i.e., with non-vanishing tangent vector) curve of class $C^k$, $k \geq 5$, in $\mathbb{RP}^2$ without inflection points. Then there exist coefficient functions $c_0,c_1,c_2$ of class $C^{k-3}$ such that $\gamma$ is the projective image of a vector-valued solution $y(t)$ of the ODE
\begin{equation} \label{general_ODE}
y'''(t) + c_2(t)y''(t) + c_1(t)y'(t) + c_0(t)y(t) = 0.
\end{equation}
By multiplying the solution $y(t)$ by a non-vanishing scalar function we may achieve that the coefficient $c_2$ vanishes identically and that $\det(y'',y',y) \equiv 1$ \cite[p.~30]{OvsienkoTabachnikov}. Subsequently decomposing the differential operator on the left-hand side of \eqref{general_ODE} in its skew-symmetric and symmetric part, we arrive at the ODE
\begin{equation} \label{special_ODE}
[y'''(t) + 2\alpha(t)y'(t) + \alpha'(t)y(t)] + \beta(t)y(t) = 0
\end{equation}
with the coefficient functions
\[ \alpha = \frac12c_1 - \frac16c_2^2 - \frac12c_2',\quad \beta = c_0 - \frac13c_1c_2 + \frac{2}{27}c_2^3 - \frac13c_2'' - \alpha'
\]
being of class $C^{k-4},C^{k-5}$, respectively \cite[p.~16]{Wilczynski}. The lift $y$ of $\gamma$ is then of class $C^{k-2}$.

The function $\beta$ transforms as the coefficient of a cubic differential $\beta(t)\,dt^3$ under reparametrizations of the curve $\gamma$. This differential is called the \emph{cubic form} of the curve \cite[pp.~15, 41]{OvsienkoTabachnikov}. Its cubic root $\sqrt[3]{\beta(t)}\,dt$ is called the \emph{projective length element}, and its integral along the curve is the \emph{projective arc length}. Points on $\gamma$ where $\beta$ vanishes are called \emph{sextactic} points. In the absence of sextactic points the curve may hence be parameterized by its projective arc length, which is equivalent to achieving $\beta \equiv 1$ and is the most natural parametrization of a curve in the projective plane \cite[p.~50]{Cartan37}.

A simple closed strictly convex curve has at least six sextactic points. This is the content of the six-vertex theorem \cite[p.~73]{OvsienkoTabachnikov} which was first proven in \cite{Mukhopadhyaya09}, according to \cite{ThorbergssonUmehara02}. Therefore such a curve does not possess a global parametrization by projective arc length.

Another common way to parameterize curves in the projective plane is the \emph{Forsyth-Laguerre} parametrization which is characterized by the condition $\alpha \equiv 0$ in \eqref{special_ODE}. This parametrization is unique up to linear-fractional transformations of the parameter $t$ \cite[pp.~25--26]{Wilczynski}, see also \cite[pp.~48--50]{Cartan37} and \cite[p.~41]{OvsienkoTabachnikov}. This implies that the curve $\gamma$ carries an invariant projective structure, which was called the \emph{projective curvature} in \cite[p.~15]{OvsienkoTabachnikov}. It is closely related to the projective curvature in the sense of \cite[p.~107]{Cartan37}, which is defined as the value of the coefficient $\alpha$ in the projective arc length parametrization. Locally projective structures on closed curves in general have been studied in \cite{Kuiper53}.

To \eqref{special_ODE} we may associate the second-order differential equation
\begin{equation} \label{reduced_ODE}
x''(t) + \frac12\alpha(t)x(t) = 0,
\end{equation}
whose solution is of class $C^{k-2}$. It is not hard to check \cite[p.~121]{Kirillov82} that if $x_1,x_2$ are linearly independent solutions of ODE \eqref{reduced_ODE}, then the products $x_1^2,x_1x_2,x_2^2$ are linearly independent $C^{k-2}$ solutions of the ODE
\begin{equation} \label{associated_ODE}
w'''(t) + 2\alpha(t)w'(t) + \alpha'(t)w(t) = 0
\end{equation}
which can be obtained from \eqref{special_ODE} by retaining the skew-symmetric part only. These solutions satisfy the homogeneous quadratic relation $x_1^2\cdot x_2^2 = (x_1x_2)^2$. Hence the vector-valued solution of ODE \eqref{associated_ODE} maps to the projective ellipse $\varsigma$ defined by this relation.

This construction is equivariant with respect to reparametrizations of the curve $\gamma$ in the following sense \cite[Theorem 1.4.3]{OvsienkoTabachnikov}.

{\lemma \label{lem:reparam_equivariance} Let $t \mapsto s(t)$ be a reparametrization of the curve $\gamma$, and let $\tilde\alpha(s)$ be the corresponding coefficient in ODE \eqref{special_ODE} in the new parameter. Let $x(t)$ be a vector-valued solution of ODE \eqref{reduced_ODE} with linearly independent components. Then there exists a non-vanishing scalar function $\sigma(s)$ such that $\tilde x(s) = \sigma(s)x(t(s))$ is a vector-valued solution of the ODE $\frac{d^2\tilde x(s)}{ds^2} + \frac12\tilde\alpha(s)\tilde x(s) = 0$. \qed}

\medskip

Obviously the scalar $\sigma(s)$ may be chosen to be positive. In fact, if we restrict to reparametrizations satisfying $\frac{ds}{dt} > 0$ and normalize the solutions $x(t),\tilde x(s)$ such that $\det(x,\frac{dx}{dt}) = \det(\tilde x,\frac{d\tilde x}{ds}) \equiv 1$, then $\sigma(s) = \sqrt{\frac{ds}{dt}}$ \cite[eq.~(2)]{LazutkinPankratova75}.

Now if two vector-valued functions $x(t),\tilde x(t)$ satisfying ODE \eqref{reduced_ODE} with coefficient functions $\alpha(t),\tilde\alpha(t)$, respectively, are related by a scalar factor, $\tilde x(t) = \sigma(t)x(t)$ for some non-vanishing $\sigma$, and $\det(x,x') = \det(\tilde x,\tilde x') \equiv 1$, then $\alpha$ and $\tilde\alpha$ coincide \cite[Theorem 1.3.1]{OvsienkoTabachnikov}. We can then reformulate above lemma as follows.

{\corollary \label{cor:reparam_equivariance} Let $\gamma(t)$ be a curve in $\mathbb{RP}^2$ without inflection points, and let $y(t)$ be a lift of $\gamma$ satisfying ODE \eqref{special_ODE} with some coefficient function $\alpha(t)$. Let $t \mapsto s(t)$ be a reparametrization of the curve $\gamma$. Let $x(t),\tilde x(s)$ be vector-valued solutions of ODE \eqref{reduced_ODE} with linearly independent components and with coefficient functions $\alpha(t),\tilde\alpha(s)$, respectively. Suppose further that $\det(x,\frac{dx}{dt}) = \det(\tilde x,\frac{d\tilde x}{ds}) \equiv 1$, and that there exists a non-vanishing scalar function $\sigma(s)$ such that $\tilde x(s) = \sigma(s)x(t(s))$ for all $s$. Then $\gamma(s)$ has a lift $\tilde y(s)$ which is a solution of ODE \eqref{special_ODE} with $\tilde\alpha(s)$ as the corresponding coefficient. \qed}

\medskip

It follows from the above that we may choose $\tilde y(s) = \sigma^2(s)y(t(s))$.

\medskip

\begin{figure}
\centering
\includegraphics[width=11.61cm,height=4.59cm]{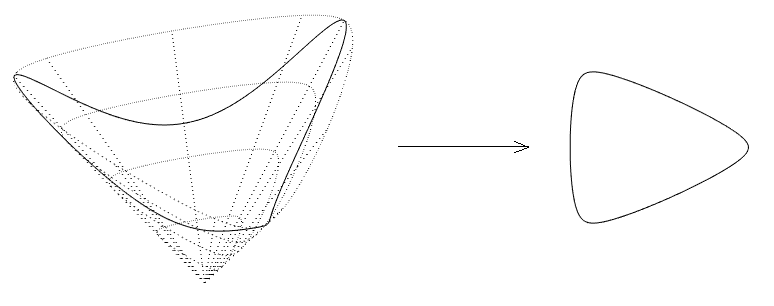}
\label{fig1}
\caption{Solution $y(t)$ on the boundary of the cone $K$ and its projection onto the simple closed convex curve $\gamma$ in an affine chart in $\mathbb RP^2$.}
\end{figure}

If $\gamma$ is represented as the projective image of a solution $y(t)$ of ODE \eqref{special_ODE}, then the dual curve $\gamma^*$ is represented as the projective image of a solution $z(t)$ of the adjoint ODE \cite[p.~61]{Wilczynski}, \cite[p.~16]{OvsienkoTabachnikov}
\begin{equation} \label{dual_ODE}
[z'''(t) + 2\alpha(t)z'(t) + \alpha'(t)z(t)] - \beta(t)z(t) = 0.
\end{equation}
Simple closed convex projective plane curves (i.e., without self-intersections, contained and convex in some affine chart on $\mathbb RP^2$) are canonically isomorphic to the manifold of boundary rays of convex proper three-dimensional cones. The solution $y(t)$ evolves on the boundary $\partial K$, while $z(t)$ evolves on $\partial K^*$, the boundary of the dual cone (see Fig.~\ref{fig1}).

Since the curve $\gamma$ is closed, we may parameterize it $2\pi$-periodically by a variable $t \in \mathbb R$. In this case the coefficient functions $\alpha(t),\beta(t)$ are also $2\pi$-periodic. The behaviour of solutions of ODEs with periodic coefficients is the subject of Floquet theory \cite{Floquet}. Namely, a shift of the variable $t$ by $2\pi$ maps the solution space of ODE \eqref{reduced_ODE} to itself, and there exists $T \in SL(2,\mathbb R)$ such that $x(t+2\pi) = Tx(t)$ for all $t \in \mathbb R$. The map $T$ is called the \emph{monodromy} of equation \eqref{reduced_ODE}. The conjugacy class of the monodromy as well as the winding number of the vector-valued solution $x(t)$ of \eqref{reduced_ODE} around the origin over one period are invariant under reparametrizations $t \mapsto s(t)$ of $\gamma$ satisfying $s(t+2\pi) = s(t)+2\pi$, i.e., preserving the periodicity condition \cite[pp.~24--25, 34--35]{OvsienkoTabachnikov}.

\medskip


Equation \eqref{reduced_ODE} with periodic coefficient function has been well studied and is known under the name \emph{Hill equation}. In \cite{LazutkinPankratova75} a complete classification of the coefficient functions under the equivalence relation generated by the group of orientation-preserving diffeomorphisms of $S^1$ and a construction of corresponding normal forms has been achieved. The equations can be classified according to several criteria. They may be divided in stable, semi-stable and unstable ones, according to the asymptotic behaviour of the solutions, or into oscillating and non-oscillating ones, according to the behaviour of the argument of the vector-valued solution. Stable solutions are always oscillating. The normal forms of the non-oscillating and the stable equations have constant coefficient functions, while in the remaining cases their coefficient functions are sinusoidal.

In \cite{Kirillov82} it has been established that the $2\pi$-periodic solutions of equation \eqref{associated_ODE} can be seen as vector fields generating diffeomorphisms of $S^1$ which preserve the coefficient function in \eqref{reduced_ODE}, and at least one non-trivial periodic solution always exists. If such a solution is nowhere zero, then it can be used to construct a diffeomorphism of $S^1$ which takes the coefficient function $\alpha$ to a constant. Moreover, this diffeomorphism is unique up to a rotation of $S^1$ if and only if $\alpha \not= \frac{n^2}{2}$ for all $n \in \mathbb N_+$. Equations with different values of the constant are non-equivalent.

Our strategy will consist in constructing diffeomorphisms of $S^1$ which transform the coefficient function of Hill equation \eqref{reduced_ODE} to a constant $\alpha \leq \frac12$. In particular, we prove the following result.

{\theorem \label{thm:existence} Let $\gamma$ be a simple closed convex projective plane curve of class $C^k$, $k \geq 5$, without inflection points. Then there exists a $2\pi$-periodic parametrization of $\gamma$ of class $C^{k-1}$ by a real variable $t$ and a $2\pi$-periodic lift $y: \mathbb R \to \mathbb R^3$ of $\gamma$ of class $C^{k-2}$ such that $y(t)$ is a solution of ODE \eqref{special_ODE} with $\alpha \equiv const$. Here the value of the constant $\alpha$ is uniquely determined by the curve $\gamma$. }

Since the classification results in \cite{LazutkinPankratova75,Kirillov82} have been established in the $C^{\infty}$ setting, we shall provide an independent proof.

\medskip

We shall now briefly summarize the contents of the paper.
First we explicitly describe the solution $z(t)$ of the adjoint ODE \eqref{dual_ODE} in terms of $y(t)$ (Lemma \ref{lem:dual}). Next we show that during each period of length $2\pi$ the projective image of the vector-valued solution $w(t)$ of ODE \eqref{associated_ODE} can make at most one turn around the ellipse $\varsigma$ on which it evolves. Equivalently, the solution $x(t)$ of ODE \eqref{reduced_ODE} can make at most one half of a turn around the origin (Lemma \ref{lem:main_obstruction}). This heavily restricts the behaviour of the solution $x(t)$ (Lemma \ref{lem:x_solution_cases}) and allows to construct a reparametrization of $\gamma$ which makes the coefficient $\alpha$ constant (Theorem \ref{thm:existence}). The value of the constant $\alpha$ depends on the eigenvalues of the monodromy $T$ of ODE \eqref{reduced_ODE} and is hence uniquely determined by the curve $\gamma$. It follows in particular that in general the Forsyth-Laguerre parametrization cannot be extended to the whole closed curve $\gamma$ (Corollary \ref{cor:no_Forsyth}).

We call a $2\pi$-periodic parametrization of $\gamma$ \emph{balanced} if the corresponding coefficient function $\alpha$ in \eqref{special_ODE} is constant.

\section{Balanced parametrizations}

Let $\gamma$ be a simple closed convex projective plane curve of class $C^k$, $k \geq 5$, without inflection points. Let the lift $y(t)$ of $\gamma$ be a $2\pi$-periodic vector-valued solution of ODE \eqref{special_ODE} such that $\det(y'',y',y) \equiv 1$. The $2\pi$-periodic coefficient functions $\alpha,\beta$ are then of class $C^{k-4},C^{k-5}$, respectively, and $y$ is of class $C^{k-2}$.

Denote $Y = (y''+\alpha y,y',y) \in SL(3,\mathbb R)$, then \eqref{special_ODE} is equivalent to the matrix-valued ODE
\begin{equation} \label{matrixODE_Y}
Y' = Y \cdot A_-,
\end{equation}
where for convenience we denoted $A_{\pm} = \begin{pmatrix} 0 & 1 & 0 \\ -\alpha & 0 & 1 \\ \pm\beta & -\alpha & 0 \end{pmatrix}$. We now describe the dual objects in terms of the matrix $Y$.

{\lemma \label{lem:dual} Assume above conditions. Let $\gamma^*$ the dual projective curve of $\gamma$. There exists a vector-valued solution $z$ of \eqref{dual_ODE} which is a lift of $\gamma^*$ and satisfies $\det(z'',z',z) \equiv 1$. The matrix $Z = (z''+\alpha z,z',z) \in SL(3,\mathbb R)$ is given by $Z = Y^{-T} Q$ with
\[ Q = \begin{pmatrix} 0 & 0 & 1 \\ 0 & -1 & 0 \\ 1 & 0 & 0 \end{pmatrix}.
\] }

\begin{proof}
Denote the matrix product $Y^{-T}Q$ by $Z$ and let $z$ be its third column. Clearly $Z$ is unimodular and $2\pi$-periodic. In particular, $z$ is non-zero everywhere. Further, by \eqref{matrixODE_Y} the product $Z$ satisfies the differential equation
\[ Z' = -Y^{-T}(YA_-)^TY^{-T}Q = -ZQ^{-1}A_-^TQ = Z \cdot A_+.
\]
It follows that $Z = (z''+\alpha z,z',z)$ and $z$ is a solution of ODE \eqref{dual_ODE}. It follows also that $\det(z'',z',z) \equiv 1$. Finally, we have $Y^TZ = Q$, which implies $\langle y(t),z(t) \rangle = \langle y'(t),z(t) \rangle = 0$ for all $t$. Hence the vector $z(t)$ is orthogonal to the plane spanned by $y(t)$ and $y'(t)$, and the projective image of $z(t)$ is the corresponding point $\gamma^*(t)$ on the dual projective curve. Thus $z$ satisfies all required conditions.
\end{proof}

Note that the dual curve $\gamma^*$ is also simple closed convex and of class $C^k$ without inflection points. Let now $t_0 \in \mathbb R$ and set $y_0 = y(t_0)$, $z_0 = z(t_0)$. Define the scalar $C^{k-2}$ functions $\mu(t) = \langle y(t),z_0 \rangle$, $\nu(t) = \langle y_0,z(t) \rangle$. By convex duality these functions are nonnegative, and $\mu(t) = 0$ or $\nu(t) = 0$ if and only if $t - t_0$ is an integer multiple of the period $2\pi$.

Assume the notations of Lemma \ref{lem:dual}. We have $ZQY^T = I$ and hence
\begin{align*}
0 &= \langle y_0,z_0 \rangle = y_0^TZQY^Tz_0 = (\nu''+\alpha\nu,\nu',\nu)Q(\mu''+\alpha\mu,\mu',\mu)^T \\ &= \nu\mu''+2\alpha\nu\mu+\mu\nu''-\nu'\mu'.
\end{align*}
For $t_0 < t < t_0+2\pi$ define the $C^{k-3}$ functions $\xi = \frac{\mu'}{\mu}$, $\theta = \frac{\nu'}{\nu}$. Dividing the above relation by $\mu\nu$ and expressing the result in terms of $\xi,\theta$ we obtain
\[ \xi' + \theta' + \xi^2 - \xi\theta + \theta^2 + 2\alpha = 0.
\]
Introducing the variable $\psi = \frac14(\xi + \theta)$ and taking into account $\xi^2 - \xi\theta + \theta^2 = 4\psi^2 + \frac34(\xi-\theta)^2$ we obtain the differential inequality
\begin{equation} \label{diff_ineq}
\psi' + \psi^2 + \frac{\alpha}{2} = -\frac{3}{16}(\xi - \theta)^2 \leq 0.
\end{equation}

{\lemma \label{lem:main_obstruction} Assume the conditions at the beginning of this section. Let $t_0 \in \mathbb R$ be arbitrary, and let $x(t)$ be a non-trivial scalar solution of ODE \eqref{reduced_ODE}. Then $x(t)$ cannot have two distinct zeros in the interval $(t_0,t_0+2\pi)$. If $x(t_0) = x(t_0+2\pi) = 0$, then $\beta \equiv 0$ and $\gamma$ is an ellipse. }

The first assertion follows by virtue of \cite[Proposition 9, p.~130]{ValleePoussin29} from the existence of a function $\psi(t)$ satisfying \eqref{diff_ineq} on $(t_0,t_0+2\pi)$. We shall, however, give an elementary proof below.

\begin{proof}
Let $t_m \in (t_0,t_0+2\pi)$ be arbitrary and define the positive function
\[ q(t) = \exp\left(\int_{t_m}^t \psi(t)\,dt\right) = \left(\frac{\mu(t)\nu(t)}{\mu(t_m)\nu(t_m)}\right)^{1/4}
\]
on $(t_0,t_0+2\pi)$, where $\psi(t)$ is the function from \eqref{diff_ineq}. Then we obtain $q'' + \frac{\alpha}{2}q = (\psi' + \psi^2 + \frac{\alpha}{2})q \leq 0$.

Let $x(t)$ be an arbitrary non-trivial solution of ODE \eqref{reduced_ODE} on $(t_0,t_0+2\pi)$ and consider the function $r(t) = x'(t)q(t) - x(t)q'(t)$. We have $r' = x''q - xq'' = -x(q'' + \frac{\alpha}{2}q)$.

Suppose for the sake of contradiction that $x(t_1) = x(t_2) = 0$ for $t_0 < t_1 < t_2 < t_0+2\pi$ and $x(t) > 0$ for all $t \in (t_1,t_2)$. Then $x'(t_1) > 0$, $x'(t_2) < 0$, and hence $r(t_1) > 0$, $r(t_2) < 0$. But $r'(t) \geq 0$ on $(t_1,t_2)$, a contradiction. The case when $x(t)$ is negative on $(t_1,t_2)$ is treated similarly. This proves the first claim.

Since $t_0$ is arbitrary, it follows that no non-trivial solution of ODE \eqref{reduced_ODE} can have two consecutive zeros at a distance strictly smaller than $2\pi$.

Let now $x(t)$ be a non-trivial solution of ODE \eqref{reduced_ODE} such that $x(t_0) = x(t_0+2\pi) = 0$. Then $x(t)$ has constant sign on $(t_0,t_0+2\pi)$, and $r'(t)$ is either nonnegative or non-positive, depending on the sign of $x$. In any case the function $r(t)$ is monotonous on $(t_0,t_0+2\pi)$. Note that $q(t)$ and $q'(t)$ can be continuously prolonged to $t_0$ and $t_0+2\pi$ and the limits of $q(t)$ vanish. We hence have $\lim_{t \to t_0}r(t) = \lim_{t \to t_0+2\pi}r(t) = 0$. It follows that $r \equiv 0$, $r' \equiv 0$, and therefore $q'' + \frac{\alpha}{2}q \equiv 0$ on $(t_0,t_0+2\pi)$. But then inequality \eqref{diff_ineq} is actually an equality and $\xi \equiv \theta$. Then there exists a constant $c > 0$ such that $\mu \equiv c\nu$.  But $\mu(t)$ is a solution of ODE \eqref{special_ODE}, while $\nu(t)$ and hence also $c\nu(t)$ is a solution of \eqref{dual_ODE}. Subtracting \eqref{dual_ODE} from \eqref{special_ODE} with $z,y$ replaced by $\mu$, respectively, we obtain $2\beta(t)\mu(t) = 0$ on $(t_0,t_0+2\pi)$. It follows that $\beta \equiv 0$, $y(t)$ is a solution of ODE \eqref{associated_ODE} and hence $\gamma$ is an ellipse. This completes the proof.
\end{proof}

Lemma \ref{lem:main_obstruction} allows to restrict the global behaviour of the solutions of ODE \eqref{reduced_ODE}.

{\lemma \label{lem:x_solution_cases} Assume the conditions at the beginning of this section. Then exactly one of the following cases holds:

\begin{itemize}
\item[(i)] There exists a solution $x(t)$ of ODE \eqref{reduced_ODE}, normalized such that $\det(x,x') \equiv 1$, that is contained in the open positive orthant and crosses each ray of this orthant exactly once, and whose monodromy equals $T = \diag(\lambda^{-1},\lambda)$ for some $\lambda > 1$.
\item[(ii)] There exists a solution $x(t)$ of ODE \eqref{reduced_ODE}, normalized such that $\det(x,x') \equiv 1$, that is contained in the open right half-plane and crosses each ray of this half-plane exactly once, and whose monodromy equals $T = \begin{pmatrix} 1 & 0 \\ 2\pi & 1 \end{pmatrix}$.
\item[(iii)] There exists a solution $x(t)$ of ODE \eqref{reduced_ODE}, normalized such that $\det(x,x') \equiv 1$, that is bounded and turns infinitely many times around the origin, and whose monodromy equals $T = \begin{pmatrix} \cos\varphi & -\sin\varphi \\ \sin\varphi & \cos\varphi \end{pmatrix}$ for some $\varphi \in (0,\pi)$. For every $t_0 \in \mathbb R$ the solution turns by an angle of $\varphi$ around the origin in the interval $[t_0,t_0+2\pi]$.
\item[(iv)] There exists a $4\pi$-periodic solution $x(t)$ of ODE \eqref{reduced_ODE}, normalized such that $\det(x,x') \equiv 1$, and whose monodromy equals $T = -I$.
\end{itemize}
The curve $\gamma$ is an ellipse if and only if case (iv) holds. }

\begin{proof}
Let $x(t)$ be an arbitrary solution of ODE \eqref{reduced_ODE} with linearly independent components, normalized such that $\det(x,x') \equiv 1$. Any other such solution can then be obtained by the action of an element of $SL(2,\mathbb R)$. The solution $x$ turns counter-clockwise around the origin and intersects every ray transversally.

\medskip

First we shall treat the case when $\gamma$ is not an ellipse. By Lemma \ref{lem:main_obstruction} every scalar solution of ODE \eqref{reduced_ODE} has its consecutive zeros placed at distances strictly larger than $2\pi$. Hence $x$ turns by an angle strictly less than $\pi$ in any time interval of length $2\pi$. In particular, it follows that the solution $x(t)$ cannot cross any 1-dimensional eigenspace of the monodromy $T$. Indeed, suppose that for some $t_0 \in \mathbb R$ the vector $x(t_0)$ is an eigenvector of $T$. Then $x(t_0+2\pi) = Tx(t_0)$ is a positive or negative multiple of $x(t_0)$, and $x$ must have made at least half of a turn around the origin in the interval $[t_0,t_0+2\pi]$, a contradiction.

We shall now distinguish several cases according to the spectrum of the monodromy $T$ of ODE \eqref{reduced_ODE}. Let $T \in SL(2,\mathbb R)$ be such that $x(t+2\pi) = Tx(t)$ for all $T$. If $\tilde x = Ax$ for some $A \in SL(2,\mathbb R)$, then $\tilde x(t+2\pi) = \tilde T\tilde x(t)$ with $\tilde T = ATA^{-1}$. We may hence conjugate $T$ with an arbitrary unimodular matrix by switching to another solution $x$.

\emph{Case 1:} The eigenvalues of $T$ are given by $\lambda,\lambda^{-1}$ for some $\lambda > 1$. By conjugation with a unimodular matrix we may achieve $T = \diag(\lambda^{-1},\lambda)$. Since $x(t)$ cannot cross the axes, it must be confined to an open quadrant. For every point $q$ in the second or fourth open quadrant the vector $Tq$ has a polar angle strictly less than that of $q$. But $x(t)$ turns in the counter-clockwise direction, and hence cannot be contained in these quadrants. By possibly multiplying $x$ by $-1$ we may hence achieve that $x$ is contained in the open positive orthant. Now for any $t_0 \in \mathbb R$ the angles of the vectors $T^kx(t_0)$ tend to $\frac{\pi}{2}$ and those of $T^{-k}x(t_0)$ to 0 as $k \to +\infty$. Therefore the angles of $x(t)$ sweep the interval $(0,\frac{\pi}{2})$ as $t$ sweeps the real line. This is the situation described in case (i) of the lemma.

\emph{Case 2:} The eigenvalues of $T$ equal 1. Since $x(t)$ cannot be an eigenvector of $T$ for any $t$, we must have $T \not= I$ and the Jordan normal form of $T$ contains a single Jordan cell. By conjugation with a unimodular matrix we may then achieve that $T = \begin{pmatrix} 1 & 0 \\ \pm2\pi & 1 \end{pmatrix}$. Since $x(t)$ cannot cross the vertical axis, it must be contained in the left or right open half-plane. By multiplying by $-1$ we may assume the solution is contained in the right half-plane. Now if the $(2,1)$ element in $T$ equals $-2\pi$, then for every point $q$ in the open right half-plane the vector $Tq$ has a polar angle strictly less than that of $q$. This is in contradiction with the counter-clockwise movement of $x$, and this case cannot appear. Hence the $(2,1)$ element in $T$ equals $2\pi$. Then for any $t_0 \in \mathbb R$ the angles of the vectors $T^kx(t_0)$ tend to $\frac{\pi}{2}$ and those of $T^{-k}x(t_0)$ to $-\frac{\pi}{2}$ as $k \to +\infty$. Therefore the angles of $x(t)$ sweep the interval $(-\frac{\pi}{2},\frac{\pi}{2})$ as $t$ sweeps the real line. This is the situation described in case (ii) of the lemma.

\emph{Case 3:} The eigenvalues of $T$ equal $e^{\pm i\varphi}$ for $\varphi \in (0,\pi)$. By conjugation with an element in $SL(2,\mathbb R)$ we may achieve that $T = \begin{pmatrix} \cos\varphi & \mp\sin\varphi \\ \pm\sin\varphi & \cos\varphi \end{pmatrix}$. If the $(2,1)$ element of $T$ has negative sign, then for every $q \not= 0$ the angle of $Tq$ equals $2\pi - \varphi$ plus the angle of $q$. Since $x$ moves counter-clockwise, it must hence sweep an angle of at least $2\pi - \varphi > \pi$ on any interval of length $2\pi$, which is not possible. Hence the $(2,1)$ element of $T$ has positive sign, and for every $q \not= 0$ the angle of $Tq$ equals $\varphi$ plus the angle of $q$. Since $x$ cannot make a complete turn around the origin in an interval of length $2\pi$, the angle swept by the solution on any such interval equals $\varphi$. Finally note that since $T$ acts by a rotation, the norm of the solution $x$ is $2\pi$-periodic and hence uniformly bounded. This is the situation described in case (iii) of the lemma.

\emph{Case 4:} The eigenvalues of $T$ equal $-1$. Similarly to Case 2 we have $T \not= -I$, and the Jordan normal form of $T$ consists of a single Jordan cell. The eigenspace to the eigenvalue $-1$ then divides $\mathbb R^2$ in two half-planes. For every $q$ in one of the open half-planes, the point $Tq$ lies in the other open half-plane. Hence the solution $x(t)$ must cross the eigenspace, leading to a contradiction. Hence this case does not occur.

\emph{Case 5:} The eigenvalues of $T$ equal $-\lambda,-\lambda^{-1}$ for some $\lambda > 1$. By conjugation with a unimodular matrix we may achieve $T = \diag(-\lambda^{-1},-\lambda)$. Similarly to Case 1 the solution $x(t)$ must then be contained in some open quadrant. But the map $T$ maps every quadrant to the opposite quadrant. Hence $x$ must cross the axes, which leads to a contradiction. Thus this case does not occur either.

\medskip

We now consider the case when $\gamma$ is an ellipse. By Lemma \ref{lem:main_obstruction} we have $\beta \equiv 0$ and \eqref{special_ODE}, \eqref{associated_ODE} represent the same ODE. Since all solutions $y$ of ODE \eqref{special_ODE} are $2\pi$-periodic, the solutions $w$ of \eqref{associated_ODE} are also $2\pi$-periodic. But the solutions $w$ are homogeneous quadratic functions of the solutions $x$ of ODE \eqref{reduced_ODE}. Hence the latter are $4\pi$-periodic, and $T^2 = I$. If $T = I$, then every two consecutive zeros of every non-trivial scalar solution of ODE \eqref{reduced_ODE} have a distance strictly smaller than $2\pi$, leading to a contradiction with Lemma \ref{lem:main_obstruction}. Hence $T = -I$, and we are in the situation described in case (iv) of the lemma.

This completes the proof.
\end{proof}

\begin{remark}
The cases i) --- iv) in the formulation of the lemma correspond to the unstable non-oscillating, semi-stable non-oscillating, stable with $\Theta < 1$, and stable with $\Theta = 1$ cases, correspondingly, in the classification in \cite{LazutkinPankratova75}. The cases 4 and 5 in the proof correspond to the semi-stable and unstable oscillating cases in \cite{LazutkinPankratova75}.
\end{remark}

{\corollary \label{cor:no_Forsyth} Assume the conditions at the beginning of this section. If the eigenvalues of the monodromy of ODE \eqref{reduced_ODE} differ from 1, then the curve $\gamma$ does not possess a global periodic Forsyth-Laguerre parametrization. }

\begin{proof}
Suppose $\gamma$ possesses a periodic Forsyth-Laguerre parametrization by a variable $s$. In this parametrization any non-zero vector-valued solution $\tilde x(s)$ of ODE \eqref{reduced_ODE} with independent components is a straight affine line, and hence sweeps a total angle of $\pi$ in the plane.

Let now $\gamma$ be parameterized $2\pi$-periodically by a variable $t$. Every non-zero vector-valued solution $x(t)$ of ODE \eqref{reduced_ODE} with independent components must also sweep a total angle of $\pi$. From Lemma \ref{lem:x_solution_cases} it follows that the monodromy of ODE \eqref{reduced_ODE} has eigenvalues equal to 1.
\end{proof}

We are now in a position to construct the reparametrization $t \mapsto s(t)$ which makes the coefficient $\alpha$ in ODE \eqref{special_ODE} constant.

\begin{proof}[of Theorem \ref{thm:existence}]
We shall begin with an arbitrary regular $2\pi$-periodic parametrization of $\gamma$ of class $C^k$. As laid out in Section \ref{sec:intro}, there exists a $2\pi$-periodic lift $y(t)$ of $\gamma$ which solves ODE \eqref{special_ODE} with some $2\pi$-periodic functions $\alpha(t)$, $\beta(t)$ of class $C^{k-4},C^{k-5}$, respectively. The coefficient function $\alpha$ gives rise to ODE \eqref{reduced_ODE}. We shall construct a $2\pi$-periodic parametrization of $\gamma$ by a new variable $s$ from the vector-valued $C^{k-2}$ solutions $x(t) = (x_1(t),x_2(t))$ of ODE \eqref{reduced_ODE} described in Lemma \ref{lem:x_solution_cases}. Note that if we write $x_1 = r\cos\phi$, $x_2 = r\sin\phi$, then the condition $\det(x,x') \equiv 1$ implies $r^2\phi' \equiv 1$ and $\phi' = r^{-1/2}$. Since $r(t)$ is of class $C^{k-2}$, the angle $\phi$ is of class $C^{k-1}$. We consider the four cases (i) --- (iv) in Lemma \ref{lem:x_solution_cases} separately.

\emph{Case (i):} Set $s(t) = \frac{\pi}{\log\lambda}\log\frac{x_2(t)}{x_1(t)}$. Note that $s$ is an analytic function of the angle $\phi$ and hence $s(t)$ is a $C^{k-1}$ function. We have $s(t+2\pi) = \frac{\pi}{\log\lambda}\log\frac{\lambda x_2(t)}{\lambda^{-1}x_1(t)} = s(t) + 2\pi$, and the new parameter $s$ parameterizes $\gamma$ $2\pi$-periodically. Set further $c = \frac{\pi}{\log\lambda} > 0$ and $\tilde\alpha = -\frac{1}{2c^2} < 0$. Then the vector-valued function $\tilde x(s) = (\sqrt{c}\lambda^{-s/2\pi},\sqrt{c}\lambda^{s/2\pi})$ obeys the differential equation $\frac{d^2\tilde x}{ds^2} + \frac{\tilde\alpha}{2}\tilde x = 0$ and we have $\frac{\tilde x_2(s(t))}{\tilde x_1(s(t))} = \lambda^{s(t)/\pi} = \frac{x_2(t)}{x_1(t)}$ for all $t$. Moreover, $\det(\tilde x,\frac{d\tilde x}{ds}) \equiv 1$. By Corollary \ref{cor:reparam_equivariance} the coefficient $\alpha$ in ODE \eqref{special_ODE} in the new coordinate $s$ identically equals the constant $\tilde\alpha$. The coefficient $\beta$ in the new variable is given by $\tilde\beta(s) = \beta(t)(\frac{ds}{dt})^{-3}$, because $\beta$ transforms as the coefficient of a cubic differential. Hence $\tilde\beta(s)$ is as $\beta(t)$ a $C^{k-5}$ function. Therefore the solution $\tilde y(s)$ of ODE \eqref{special_ODE} in the variable $s$ is of class $C^{k-2}$.

\emph{Case (ii):} Set $s(t) = \frac{x(t_2)}{x(t_1)}$. Again $s$ is an analytic function of the angle $\phi$ and $s(t)$ is a $C^{k-1}$ function. We have $s(t+2\pi) = \frac{2\pi x(t_1) + x(t_2)}{x(t_1)} = s(t) + 2\pi$, and $s$ parameterizes $\gamma$ $2\pi$-periodically. Define $\tilde x(s) = (1,s)$, then $\det(\tilde x,\frac{d\tilde x}{ds}) \equiv 1$, $\frac{d^2\tilde x}{ds^2} = 0$, and $\frac{\tilde x_2(s)}{\tilde x_1(s)} = \frac{x(t_2)}{x(t_1)}$. By Corollary \ref{cor:reparam_equivariance} the coefficient $\alpha$ in ODE \eqref{special_ODE} in the new coordinate $s$ identically equals zero. As in the previous case the coefficient $\tilde\beta(s)$ is a $C^{k-5}$ function and the solution $\tilde y(s)$ of ODE \eqref{special_ODE} in the variable $s$ is of class $C^{k-2}$.

\emph{Case (iii):} Set $s(t) = \frac{2\pi}{\varphi}\phi(t)$. Again $s$ is a $C^{k-1}$ function and $s(t+2\pi) = \frac{2\pi}{\varphi}(\phi(t) + \varphi) = s(t) + 2\pi$, and $s$ parameterizes $\gamma$ $2\pi$-periodically. Define $c = \frac{2\pi}{\varphi}$, $\tilde\alpha = \frac{2}{c^2}$, and $\tilde x(s) = (\sqrt{c}\cos\frac{s}{c},\sqrt{c}\sin\frac{s}{c})$. Then $\det(\tilde x,\frac{d\tilde x}{ds}) \equiv 1$, $\frac{d^2\tilde x}{ds^2} + \frac{\tilde\alpha}{2}\tilde x = 0$, and the angles of $x(t)$ and $\tilde x(s)$ both equal $\phi$. By Corollary \ref{cor:reparam_equivariance} the coefficient $\alpha$ in ODE \eqref{special_ODE} in the new coordinate $s$ identically equals the constant $\tilde\alpha$. As in the previous case the coefficient $\tilde\beta(s)$ is a $C^{k-5}$ function and the solution $\tilde y(s)$ of ODE \eqref{special_ODE} in the variable $s$ is of class $C^{k-2}$.

\emph{Case (iv):} The curve $\gamma$ is an ellipse, and by an appropriate choice of the coordinate basis in $\mathbb R^3$ we may achieve that $\gamma$ is the projective image of the vector-valued function $y(t) = (1,\cos t,\sin t)$. This function is a solution of ODE \eqref{special_ODE} with $\alpha \equiv \frac12$, $\beta \equiv 0$, and the variable $t$ parameterizes $\gamma$ analytically and $2\pi$-periodically.

\medskip

Finally we show that the value of the constant $\alpha$ is uniquely determined by $\gamma$. Let the lift $y(t)$ of $\gamma$ be a $2\pi$-periodic solution of ODE \eqref{special_ODE} with constant coefficient $\alpha$. Let $x(t)$ be the solution from Lemma \ref{lem:x_solution_cases}.

If $\alpha < 0$, then $x(t)$ must be a hyperbola, hence case (i) is realized, and $\alpha$ relates to the spectrum of the monodromy $T$ of ODE \eqref{reduced_ODE} by $\alpha = -\frac{\log^2\lambda}{2\pi^2}$.

If $\alpha = 0$, then by Corollary \ref{cor:no_Forsyth} the eigenvalues of $T$ equal 1.

If $\alpha \in (0,\frac12)$, then $x(t)$ must be an ellipse and sweeps an angle strictly less than $\pi$ in any interval of length $2\pi$. Hence case (iii) is realized, and $\alpha$ is related to the spectrum of $T$ by $\alpha = \frac{\varphi^2}{2\pi^2}$.

If $\alpha \geq \frac12$, then $x(t)$ must also be an ellipse and sweeps an angle of at least $\pi$ in any interval of length $2\pi$. Hence case (iv) is realized, $x(t)$ sweeps an angle of exactly $\pi$, and $\alpha = \frac12$.

In any case $\alpha$ is uniquely determined by the spectrum of $T$. However, the spectrum of $T$ depends only on $\gamma$. Hence $\alpha$ is also uniquely determined by $\gamma$.
\end{proof}

{\definition Let $\gamma$ be a simple closed convex projective plane curve of class $C^k$, $k \geq 5$, without inflection points. We call a $2\pi$-periodic parametrization of $\gamma$ by a real variable $t$ \emph{balanced} if there exists a $2\pi$-periodic lift $y(t)$ of $\gamma$ to $\mathbb R^3$ which is a vector-valued solution of ODE \eqref{special_ODE} with $\alpha \equiv const$. }

\medskip

By Theorem \ref{thm:existence} a balanced parametrization always exists. In the case of non-quadratic curves the balanced parametrization is unique up to a shift of the variable $t$ by \cite[Lemma 2]{Kirillov82}, and hence defines an invariant metric on the curve.
For an ellipse every two balanced parametrizations are related by a projective transformation.

\section*{Acknowledgements}
The author would like to thank the referee for a thorough review and for pointing out relevant literature on normal forms of ODEs.

\bibliographystyle{plain}
\bibliography{geometry,diff_eq}

\end{document}